\documentclass[3p]{elsarticle}

\usepackage{lineno,hyperref}
\usepackage{amsmath}
\usepackage{subcaption}
\usepackage{graphicx} 
\modulolinenumbers[5]
\usepackage{xcolor}
\usepackage{soul}

\begin{document}

\begin{frontmatter}

  \title{ Stable 3D FDTD Method for arbitrary Fully Electric and Magnetic Anisotropic  Maxwell Equations}

\author[ACMS,COS]{J. Nehls \corref{mycorrespondingauthor}}
\cortext[mycorrespondingauthor]{Corresponding author}
\ead{nehls@email.arizona.edu}
\author[ACMS,DoM]{C. Dineen}
\author[NC]{J. Liu}
\author[WI]{C. Poole}
\author[ACMS,DoM]{M. Brio}
\author[ACMS,COS,DoM]{J.V. Moloney}

\address[ACMS]{Arizona Center for Mathematical Sciences, University of Arizona, Tucson, Arizona 85721}
\address[DoM]{Department of Mathematics, University of Arizona, Tucson, Arizona 85721 }
\address[COS]{College of Optical Sciences, University of Arizona, Tucson, Arizona 85721 }
\address[NC]{Department of Mathematical Sciences, Delaware State University, Dover, DE 19901} 
\address[WI]{Department of Physics, Lawrence University, Appleton, WI 54911}




\begin{abstract}
We have developed a new fully anisotropic 3D FDTD Maxwell solver for arbitrary electrically and magnetically anisotropic media for piecewise constant electric and magnetic materials that are co-located over the primary computational cells. Two numerical methods were developed that are called non-averaged and averaged methods, respectively. The non-averaged method is first order accurate, while the averaged method is second order accurate for smoothly-varying materials and reduces to first order for discontinuous material distributions. For the standard FDTD field locations with the co-location of the electric and magnetic materials at the primary computational cells, the averaged method require development of the different inversion algorithms of the constitutive relations for the electric and magnetic fields.  We provide a mathematically rigorous stability proof followed by extensive numerical testing that includes long-time integration, eigenvalue analysis, tests with extreme, randomly placed material parameters, and various boundary conditions. For accuracy evaluation we have constructed a test case with an explicit analytic solution. Using transformation optics, we have constructed complex, spatially inhomogeneous geometrical object with fully anisotropic materials and a large dynamic range of $\underline{\epsilon}$ and $\underline{\mu}$, such that a plane wave incident on the object is perfectly reconstructed downstream. In our implementation, the considerable increase in accuracy of the averaged method only increases the computational run time by 20\%.

\end{abstract}

\begin{keyword}
3D FDTD \sep Fully-Anisotropic\sep Stability \sep Eigenvalue Analysis \sep Cloaking
\MSC[2010] 00-01\sep  99-00
\end{keyword}

\end{frontmatter}

\section{Introduction}
\label{sec:intro}
A rapidly developing field of metamaterials, multiferroic and magnetoelectric materials, applications of transformation optics are at the forefront of the important advancements, new discoveries and practical applications in electomagnetism, \cite{Meta}, \cite{magnetophotonics}, \cite{TOptics}, \cite{Benson}, \cite{magnonic}. Computer simulations of such materials require robust and accurate numerical Maxwell solvers when both electric and magnetic materials are fully anisotropic materials.

The finite-difference time-domain (FDTD) method has a long history of success with isotropic and diagonally-anisotropic materials \cite{Taflove05}. 
The key advantages of FDTD are its second-order accuracy, enforcement of continuity of the respective normal and tangential field components across interfaces on Cartesian grids, and its non-dissipative energy preserving nature.

The previous work primarily concentrated on electrically anisotropic materials, \cite{choi1986finite}, \cite{schneider1993finite}, \cite{zhao1999efficient},   \cite{werner2007stable}  and \cite{werner2013more}. A split-step FDTD method for 3D Maxwell’s equations for the fully anisotropic, but homogeneous media, is presented in Singh \emph{et al.} \cite{singh2010split}. Our work is inspired by the  \cite{werner2007stable}  and \cite{werner2013more}. In the mathematical proof, we have utilized their idea of "triplets", but  without introducing new terminology and more importantly, the proof switches from the domain of dependence point of view to the domain of influence argument as explained in section \ref{sec:SPD}.  We have generalizes the results of these two papers on FDTD that deal with anisotropic permittivity only and semi-discrete stability considerations to fully electric and magnetic anisotropic materials and provided a rigorous mathematical proof of stability in a fully discrete case. We have also shown that stability of a particular algorithm may depend on the material placement. In particular, "stable" algorithm of \cite{werner2007stable} that was determined to be unstable and modified to be stable in  \cite{werner2013more}, are all the same algorithm in our material placement when applied to the electric field, and is proven to be stable in fully discrete case. The magnetic field treatment requires a different algorithm to be provable stable. We have also added a missing component in the proof in  \cite{werner2013more} that the global material matrix is an SPD matrix by including necessary permutations applied to global electric field and electric displacement vectors, as well as respective material matrices. 

 In this paper we have developed new fully anisotropic 3D FDTD Maxwell solver for arbitrary electrically and magnetically anisotropic media for piecewise constant electric and magnetic materials that are co-located at the primary computational cells. In particular, we provide a mathematically rigorous proof of the stability of the method. Extensive numerical testing involving extreme and random parameter regimes, and various boundary conditions are utilized to verify the correctness of the implementation and to illustrate the performance of the method. For accuracy evaluation we have constructed a test case with an explicit analytic solution. Using transformation optics, we have constructed complex, spatially inhomogeneous geometrical object with fully anisotropic materials and a large dynamic range of $\underline{\epsilon}$ and $\underline{\mu}$, such that a plane wave incident on the object is perfectly reconstructed downstream.
 
  The co-location of electric and magnetic materials at the primary computational cells is done for the following reasons: a) since physical fully anisotropic objects do not have electric and magnetic parts shifted by half of the computational cell this approach avoids the necessity for complicated cut-cell algorithms;  b) to be suitable for generalizations and applications in finite element context utilizing conventional mesh generators where material discontinuities are allowed across the boundaries between the elements, but not within the individual elements itself. Note, that the co-location of the materials over the primary computational cells breaks the symmetry between the electric and magnetic materials/fields placement. Therefore, two separate algorithms for inverting the constitutive relations were developed, one for the electric field and another one for the magnetic field. 

The paper is organized as follows: in section \ref{sec:aFDTD}, we present the fully anisotropic FDTD algorithm. 
In section \ref{sec:astability}, we analyze the stability of the anisotropic FDTD algorithm, followed by the 
numerical stability tests in section \ref{sec:nstability}. The accuracy study is presented in section \ref{sec:accuracy}. 

\section{Fully anisotropic FDTD algorithm}
\label{sec:aFDTD}
\begin{figure}[htb] 
\includegraphics[width=100px]{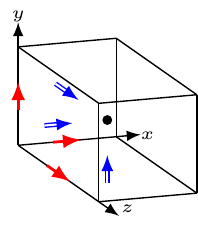}
\centering\caption{Location of the $E$-field components on the faces (double blue arrows) and the $H$-field components on the edges (red arrows) and piecewise materials ($\underline{\epsilon} \text{ and } \underline{\mu}$) location marked at cell center (black dot). The diagram is an example of a single computational cell where each field and material shown is indexed in $x, y,$ and $z$ directions using discrete indices, $i, j,$ and $k$, respectively.}
\label{fig:yeeCube}
\end{figure}

In our implementation of the 3D FDTD method, the relationship between material locations and field components is as follows: for a given computational cell with the origin at the vertex  $(i,j,k)$, the $D/E$-field components are located at the center of their respective low-side faces, while the $B/H$-field components are located on the low-side of their respective edges, and both the permittivity and permeability matrices are assumed be piecewise constant over the primary computational cell, as shown in figure \ref{fig:yeeCube}. For example, The $D_{x(i,j,k)}$ and $E_{x(i,j,k)}$
 components attributed to computational cell $(i,j,k)$ are located at $(x_{i}, y_{j+1/2}, z_{k+1/2})$, while  $B_{x(i,j,k)}$ and $H_{x(i,j,k)}$   are located at $(x_{i+1/2}, y_{j}, z_{k})$. 
  
If the material distribution is described analytically, either by a smooth or non-smooth function, the material values are replaced by their average values over the primary computational cells. 

The standard FDTD discretization of Maxwell's equations is as follows:

\begin{eqnarray} \label{eqn:yee}
\textbf{B}^{n+1/2} &=& \textbf{B}^{n-1/2} - \Delta t \: \nabla \times (\textbf{E}^{n}), \\
\textbf{H}^{n+1/2} &=& \underline{\zeta}\; \textbf{B}^{n+1/2},  \label{eq:const1} \\
\textbf{D}^{n+1} &=& \textbf{D}^{n} +  \Delta t \: \nabla \times (\textbf{H}^{n+1/2}), \\
\textbf{E}^{n+1} &=& \underline{\xi}\; \textbf{D}^{n+1}, \label{eq:const2}
\end{eqnarray} 
where $\underline{\zeta}$ and $\underline{\xi}$ are the inverses of the   $3\times 3$ symmetric positive definite matrices $\underline{\mu}$ and $\underline{\epsilon}$, respectively. This method has been shown to be stable for the diagonally-anisotropic, homogeneous  $\underline{\epsilon}$ and $\underline{\mu}$, \cite{remis2000stability}.
For fully anisotropic materials, the only changes to the standard FDTD method are due to the presence of the fully anisotropic material constitutive relations, equations \eqref{eq:const1} and \eqref{eq:const2}.   
 
\subsection{Non-Averaged Fully Anisotropic Method} \label{sec:nonAve}
The constitutive relations introduce a first-order error as non co-located field components assigned to a single computational cell ($i,j,k$) are used, and they are on the order of cell size apart from each other.

For example the $E_z$ update for the computational cell $(i, j, k)$ is 

\begin{equation}
  E_{z (i,j,k)} =  \left[ \begin{array}{c}  \xi_{zx} \enskip  \xi_{zy} \enskip  \xi_{zz} \end{array} \right]_{ijk} 
                  \left[ \begin{array}{c}  D_x \\     D_y  \\  D_z \end{array} \right]_{ijk}. 
\end{equation}

The  update for all other field components is shown  in \ref{ap:non-centeredUpdate}.

\subsection{Averaged Fully Anisotropic Method}
Here we define two separate update methods, one for the $E$ and one for the $H$ field components.  
This is due to the difference in the location of the $E$ and the $H$ fields that are located  at the primary computational cell faces and the computational cell edges, respectively, while both  materials ($\underline{\epsilon} \text{ and } \underline{\mu}$)  are co-located at the primary computational cells as illustrated in figure \ref{fig:yeeCube}.
While being different in terms of number of surrounding computational cells, both update methods follow a similar strategy: 1) first, average $D$s/$B$s inside the computational cells that share a common face (edge) of interest; 2) invert respective constitutive relations in each of the surrounding computational cells; 3) take the arithmetic average of the resulting $E$/$H$ field components. 
In the following subsections, we describe both algorithms for $E_z$ and $H_z$ components. The update algorithms for all of the $E$ and $H$ field components are given  in \ref{ap:aveUpdateE} and \ref{ap:aveUpdateH}.

\subsubsection{E-Update}
\begin{figure}[htb] 
\includegraphics[width=100px]{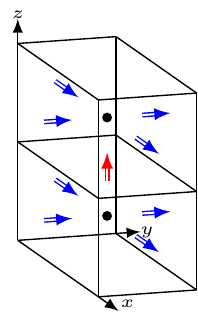}
\centering\caption{ The field components and materials associated  $E_z/D_z$  component (red).}
 
\label{fig:Etriplets}
\end{figure}

In order to update $E_{z(i,j,k)}$, consider two computational cells that share common face where $E_{z(i,j,k)}$ is located, as shown in figure \ref{fig:Etriplets}.

As mentioned above, computing first averaged values of $D$s in each computational cell, followed by inverting the constitutive relations in each computational cell, and finally, taking the arithmetic average of the resulting $E_z$ values, gives the following update formula for $E_{z(i,j,k)}$,

\begin{small}
\begin{eqnarray}\label{eq:EzAveUpdate}
  E_{z(i,j, k)} &=& \frac{1}{8}  \left[ \begin{array}{c}  \xi_{zx} \enskip  \xi_{zy} \enskip  \xi_{zz} \end{array} \right]_{i,j,k} 
                       \left(     \left[ \begin{array}{c}  D_{x(i,j,k)} \\     D_{y(i,j,k)}  \\  D_{z(i,j,k)} \end{array} \right]   +
  \left[ \begin{array}{c} D_{x(i+1,j,k)} \\ D_{y(i,j,k)}  \\  D_{z(i,j,k)} \end{array} \right] +
  \left[ \begin{array}{c} D_{x(i,j,k)} \\   D_{y(i,j+1,k) } \\ D_{z(i,j,k)}   \end{array} \right]  +
  \left[ \begin{array}{c} D_{x(i+1,j,k)} \\ D_{y(i,j+1,k) } \\ D_{z (i,j,k) }   \end{array} \right]   \right) \\
                &+& \frac{1}{8} \left[ \begin{array}{c}  \xi_{zx} \enskip  \xi_{zy} \enskip  \xi_{zz} \end{array} \right]_{ij,k-1}
                   \left(     \left[ \begin{array}{c}  D_{x(i,j,k-1)} \\     D_{y(i,j,k-1)}  \\  D_{z(i,j,k)} \end{array} \right]   +
  \left[ \begin{array}{c} D_{x(i+1, j,k-1)} \\ D_{y(i,j,k-1)}  \\  D_{z(i,j,k)} \end{array} \right] +
  \left[ \begin{array}{c} D_{x(i,j,k-1)} \\   D_{y(i,j+1,k-1) } \\ D_{z(i,j,k)}   \end{array} \right]  +
  \left[ \begin{array}{c} D_{x(i+1,j,k-1)} \\ D_{y(i,j+1,k-1) } \\ D_{z(i,j,k) }   \end{array} \right]   \right).   \nonumber
\end{eqnarray}
\end{small}
For implementation efficiency the common terms are grouped together,  
\begin{eqnarray}\label{eq:ezImp}
E_{z(i,j,k)} &=& \frac{1}{4} \left( \xi_{zx(i,j,k)} (D_{x(i,j,k)} + D_{x(i+1,j,k)}) +        \xi_{zx(i,j,k-1)} (D_{x(i,j,k-1)} + D_{x(i+1,k,k-1)})\right)\\ \nonumber
         &+& \frac{1}{4}\left( \xi_{zy(i,j,k)} (D_{y(i,j,k)} + D_{y(i,j+1,k)}) +    \xi_{zy(i,j,k-1)} (D_{y(i,j,k-1)} + D_{y(i,j+1,k-1)})\right)\\ \nonumber
           &+& \frac{1}{2}\left(\xi_{zz(i,j,k)} + \xi_{zz(i,j,k-1)}\right) D_{z(i,j,k)}.\\ \nonumber
\end{eqnarray}

\subsubsection{H-Update}
To update the edge-centered $H_{z(i,j,k)}$ field component, consider the four computational cells that share common edge where $H_{z(i,j,k)}$ is located with materials $\zeta_{(i,j,k)}$, $\zeta_{(i-1,j,k)}$,  $\zeta_{(i,j-1,k)}$, and  $\zeta_{(i-1,j-1,k)}$, respectively. 

\begin{figure}[htb] 
\includegraphics[width=200px]{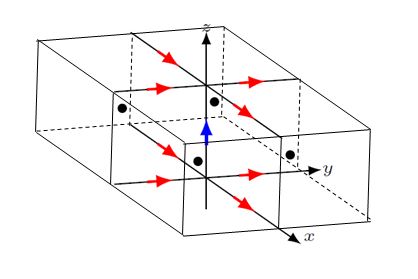}
\centering\caption{  The field components and materials associated  $H_z/B_z$  component (blue).}
\label{fig:Htriplets}
\end{figure}

Since there are four surrounding computational cells for each magnetic field component (vs two for each electric field component), we compute first the averaged values of $B$s in each of the four surrounding computational cells. This is followed by inverting the constitutive relations within each computational cell, and finally, taking the arithmetic average of the resulting $H_z$ values,  gives the following update formula for $H_{z(i,j,k)}$,

\begin{small}
\begin{eqnarray}\label{eq:HzAveUpdate}
  H_{z(i,j, k)} &=& \frac{1}{8}  \left[ \begin{array}{c}  \zeta_{zx} \enskip  \zeta_{zy} \enskip  \zeta_{zz} \end{array} \right]_{i,j,k} 
                   \left(     \left[ \begin{array}{c}  B_{x(i,j,k)} \\     B_{y(i,j,k)}  \\  B_{z(i,j,k)} \end{array} \right]   +
  \left[ \begin{array}{c} B_{x(i, j,k+1)} \\ B_{y(i,j,k+1)}  \\  B_{z(i,j,k)} \end{array} \right] \right) \\ \nonumber
               &+&  \frac{1}{8} \left[ \begin{array}{c}  \zeta_{zx} \enskip  \zeta_{zy} \enskip  \zeta_{zz} \end{array} \right]_{i-1,j,k}
                   \left(     \left[ \begin{array}{c} B_{x(i-1,j,k)} \\   B_{y(i,j,k) } \\ B_{z(i,j,k)}   \end{array} \right]  +
                      \left[ \begin{array}{c} B_{x(i-1,j,k+1)} \\   B_{y(i,j,k+1) } \\ B_{z(i,j,k)}   \end{array} \right] \right) \\ \nonumber
 &+& \frac{1}{8}  \left[ \begin{array}{c}  \zeta_{zx} \enskip  \zeta_{zy} \enskip  \zeta_{zz} \end{array} \right]_{i,j-1,k} 
                   \left(     \left[ \begin{array}{c}  B_{x(i,jk)} \\     B_{y(i,j-1,k)}  \\  B_{z(i,jk)} \end{array} \right]   +
  \left[ \begin{array}{c} B_{x(i, j,k +1)} \\ B_{y(i,j-1,k+1)}  \\  B_{z(i,jk)} \end{array} \right] \right) \\ \nonumber
               &+&  \frac{1}{8} \left[ \begin{array}{c}  \zeta_{zx} \enskip  \zeta_{zy} \enskip  \zeta_{zz} \end{array} \right]_{i-1,j-1,k}
                   \left( \left[ \begin{array}{c} B_{x(i-1,j,k)} \\   B_{y(i,j-1,k) } \\ B_{z(i,j,k)}   \end{array} \right]  +
  \left[ \begin{array}{c} B_{x(i-1,j,k+1)} \\ B_{y(i,j-1,k+1) } \\ B_{z(i,jk)}   \end{array} \right]  \right). \\  \nonumber
\end{eqnarray}
\end{small}
The simplification of this expression as well as the update algorithms for all of the $H$ field components are given   in \ref{ap:aveUpdateH}.

\section{Stability Analysis}
\label{sec:astability}

In this section, we prove the stability of the two fully anisotropic FDTD methods described in the previous section. We apply a combined approach of Fourier harmonic ansatz (von Neumann analysis) in time and a matrix stability analysis in space. In \cite{GedneyRoden}, under an additional assumption that the curl-curl matrix is diagonalizable, it was shown that requiring that the material matrices are symmetric and positive definite (SPD) is sufficient for the stability under the appropriate CFL restriction.  Below, we prove that both methods of the previous section have a diagonalizable curl-curl matrix as well as SPD material matrices. Therefore, similarly to the standard FDTD algorithm, both methods are neutrally stable, e.g. $\omega$ is real for each Fourier mode, and thus there is neither growth nor decay in time, under the CFL restriction.

Consider a fully anisotropic Maxwell's equations in non-dispersive media ($J=0$), 
\begin{eqnarray} \label{eqn:constEqn}
  D_t &=& \nabla \times (M_{\zeta} B),  \\ 
B_t &=& - \nabla \times  (M_{\xi} D),  
\end{eqnarray}
where $B$ and $D$ are electric and magnetic fluxes, and $M_{\zeta}$ and $M_{\xi}$  are the inverses of the permittivity and permeability matrices, respectively,
\begin{eqnarray}\label{eq:constE}
E &=& M_{\xi} D,   \\ 
H &=& M_{\zeta} B. \label{eq:constH}
\end{eqnarray}
Eliminating $B$ from the above equations gives 
$ D_{tt} = -\nabla \times (M_{\zeta} \nabla \times (M_{\xi} D))$. Assuming a harmonic ansatz in time,
$\displaystyle D(x,y,z,t) = e^{ i \omega t } \tilde D(x,y,z)$, results in the following eigenvalue problem 
\begin{equation}
\omega^2 \tilde D = \nabla \times (M_{\zeta} \nabla \times (M_{\xi} \tilde D)). 
\end{equation}
To eliminate exponentially growing solutions $\omega $ should be real, and therefore, the eigenvalues of the right-hand side curl-curl operator should be non-negative. For this it is sufficient to require that the material matrices are SPD matrices, \cite{Cohen2017}.  

A similar argument is applied to the fully discrete case below. 
Consider an FDTD algorithm,  
\begin{eqnarray}
\frac {D^n-D^{n-1}}{\Delta t}  &=& C_h M_{\zeta} B^{n -{1/2}}, \label{eqn:Dupdate} \\
\frac {B^{n+1/2}-B^{n-1/2}}{\Delta t} &=& -  C_e M_{\xi} D^n, \label{eqn:Bupdate}
\end{eqnarray}
where $C_h$ and $C_e$ are the discrete curl operators applied to the electric and magnetic fields, respectively. 
Note that on a uniform grid, for PEC or periodic boundary conditions, the following reciprocity relation holds, $C_h=C_e ^T$ \cite{weiland1984numerical}.    

To eliminate $B$, we take the forward time difference of equation (\ref{eqn:Dupdate}) and substitute  $(B^{n+1/2}-B^{n-1/2} )/\Delta t$ using equation \ref{eqn:Bupdate}. This gives, 
\begin{equation}
(D^{n+1}- 2 D_n +D^{n-1} )/\Delta t^2   = - C_h M_{\zeta} C_e M_{\xi} D^n. 
\end{equation}
Finally, assuming periodicity in time and applying Fourier harmonic ansatz, $\displaystyle D^n=\tilde D \hspace{1mm} e^{ i \omega n \Delta t }$, results in the following eigenvalue problem, 
\begin{equation}
  sin^2 (\omega \Delta t /2)   \tilde D =  \frac { \Delta t^2} {4}   C_h M_{\zeta} C_e M_{\xi}  \tilde D. \label{eqn:Disper}
\end{equation}
As for the continuous case, for the stability requirement on $\omega $ to be real, it is suffices that the right-hand side (discrete curl-curl operator) is diagonalizable and has positive eigenvalues, under the CFL time step restriction that is expressed in terms of the spectral radius of the discrete curl-curl operator in \cite{GedneyRoden, remis2000stability, denecker2004new},
\begin{equation}
\Delta t \le \frac {2}{ \sqrt {\rho (C_h M_{\zeta} C_e M_{\xi}) }},
\label{eqn:CFLequ}
\end{equation}
where $\rho$ is the spectral radius of the matrix $C_h M_{\zeta} C_e M_{\xi}$.
To show that SPD material matrices imply that the discrete curl-curl operator is diagonalizable and has non-negative eigenvalues (assuming,  the reciprocity relation $C_e = C, C_h = C^T$ applies), first multiply both sides of equation (\ref{eqn:Disper}) by   
$M_{\xi}^{1/2}$. Then using the fact that the material matrices are SPD,  $\displaystyle M_{\xi}=M_{\xi}^{1/2} M_{\xi} ^{1/2}$ , where the square root matrices themselves are SPD matrices \cite{werner2013more, weiland1984numerical}, equation (\ref{eqn:Disper}) can be written as 
\begin{equation}
sin^2 (\frac {\omega \Delta t}{2})  M_{\xi}^{1/2}  \tilde D = \frac { \Delta t^2} {4}   (M_{\xi}^{1/2})^ T  C^T (M_{\zeta}^{1/2})^T    M_{\zeta}^{1/2} C M_{\xi}^{1/2} M_{\xi} ^{1/2}   \tilde D. 
\end{equation}
Finally, introducing $A = M_{\zeta}^{1/2} C M_{\xi}^{1/2}$ and 
$\hat D =  M_{\xi}^{1/2}  \tilde D$,  gives 
\begin{equation}
sin^2 (\frac {\omega \Delta t}{2})\; \hat D  = (A^T A)\; \hat D.  
\end{equation}
Since $A^T A$ is symmetric positive semi-definite matrix, it is diagonalizable and has non-negative eigenvalues. In the remaining part of this section we prove that both of our methods have SPD material matrices.

\subsection{SPD of the Material Matrix for the Non-Averaged Method}

 For the first order non-averaged method,  the global  material matrix is block diagonal with each block consisting of physical, $3\times 3$, SPD electric permittivity or magnetic permeability matrices with blocks ordered according to the computational cell ordering  $(i,j,k)$, where the indices are traversing along $x$, then $y$, and finally along the $z$ directions, respectively, as illustrated in equation \ref{eq:zeroUpdate} below, 


\begin{equation}\label{eq:zeroUpdate}
   \left[ \begin{array}{c}  E_{(i,j,k)} \\ E_{(i+1,j,k)} \\ \vdots \\  E_{(nx, ny, nz)}  \end{array} \right]
    = \begin{bmatrix} 
      \xi_{(i,j,k)} & 0 & \dots & 0  \\
      0 &  \xi_{(i+1,j,k)} & \dots & 0  \\
       \vdots & \vdots & \ddots & \vdots  \\
      0 & 0 & \dots & \xi_{(nx,ny,nz)}  
\end{bmatrix}
 \left[ \begin{array}{c}  D_{(i,j,k)} \\ D_{(i+1,j,k)} \\ \vdots \\  D_{(nx, ny, nz)}  \end{array} \right]
\end{equation}
where  $E_{(i,j,k)}$, $D_{(ijk)}$, and $\xi_(i,j,k)$ denote the electric field components and the inverse permittivity matrix associated with the computational cell $(i,j,k)$.  For example, 

\begin{equation}\label{eq:NonAveUpdate}
 D_{(ijk)}=\left[ \begin{array}{c}  D_x \\     D_y  \\  D_z \end{array} \right]_{ijk}, \;\;\;\;\;\; \xi_{(ijk)}= \begin{bmatrix} 
           \xi_{xx} & \xi_{xy} & \xi_{xz}\\
           \xi_{yx} & \xi_{yy} & \xi_{yz} \\
           \xi_{zx} & \xi_{zy} & \xi_{zz}
           \end{bmatrix}_{i,j,k}.
\end{equation}
The indices $(nx, ny, nz)$ label the last computational cell.

\subsection{SPD of the Material Matrix for the Averaged Method}
\label{sec:SPD}

Consider constitutive relation for electric field,  equation \eqref{eq:constE},
\begin{equation} \label{eq:Meps}
  E = M_{\xi}  D, 
\end{equation}  
where $ E$ and $ D$ are column vectors containing all field components of the 3D computational domain written in some fixed  "standard" order. For example, we may chose the standard ordering that corresponds to the computational cell $(i,j,k)$ by choosing the electric field components located on  the low-side faces of each computational cell as described in section \ref{sec:aFDTD}, see figure \ref{fig:yeeCube} and equation \eqref{eq:zeroUpdate}.  

\begin{figure}[htb!]
\begin{subfigure}{0.4\textwidth}
\includegraphics[width=0.7\linewidth]{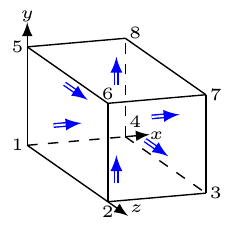}
\end{subfigure}
\begin{subfigure}{0.4\textwidth}
\includegraphics[width=0.7\linewidth]{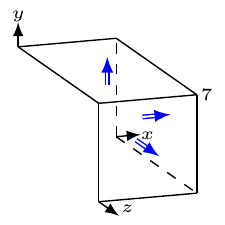}
\end{subfigure}
\caption{Enumeration of the vertices for each computational cell (Left) and the location of the electric field components at the centers of the faces corresponding to vertex 7 (Right).}
\label{fig:Vertex}
\end{figure}

Now, consider the local vertex labeling for each computational cell as indicated in figure  \ref{fig:Vertex}. Note, that while we chose to associate with each computational cell the fields located on the lower-side faces, corresponding to vertex 1 in figure  \ref{fig:Vertex}, we could have chosen any vertex and respective faces to represent the electric field of the computational cell. We enumerate these eight possibilities to represent the global electric field as  $\widetilde{E}^{(m)}$, and $\widetilde{D}^{(m)}\;m=1,2,\dots , 8$. These vectors will contain all the field components without repetition. When the standard ordering of the global electric field for each computational cell is based on vertex 1 , the global vector $ E$ in equation \eqref{eq:Meps} is the same as  $\widetilde{E}^{(1)}$. 

Equation \eqref{eq:EzAveUpdate} (and similar equations for the other field components) may be interpreted either from the "domain of dependence" point of view, where we focus on the contributions to each field component from the surrounding displacement vectors that share common face, or from the "domain of influence" argument, where we collect the electric field components that are influenced by a single displacement vector corresponding to a particular vertex within each computational cell. In fact, for any vertex of the computational cell, the corresponding displacement vector will only contribute to the electric field components that are located at the same faces. 

Consider an arbitrary computational cell and eight local displacement vectors corresponding to each vertex. The components of these vectors are located on the cell faces intersecting at a particular vertex as shown in figure \ref{fig:Vertex} for vertex 7. 


The displacement vector corresponding to vertex 7 contributes to the components of the global electric field that are located at the same faces,  

\begin{equation}
  \left[ \begin{array}{c} E_{x} \\ E_{y} \\ E_{z} \end{array} \right]_{i+1,j+1,k+1}  \Longleftarrow   \begin{bmatrix} 
          \xi_{xx} & \xi_{xy} & \xi_{xz}\\
          \xi_{yx} & \xi_{yy} & \xi_{yz} \\
          \xi_{zx} & \xi_{zy} & \xi_{zz}
          \end{bmatrix}_{i,j,k} 
                  \left[ \begin{array}{c} D_{x} \\ D_{y} \\ D_{z} \end{array} \right]_{i+1,j+1,k+1}. 
\end{equation}

Combining these relations for each computational cell into global vectors  $\widetilde{E}^{(7)}$ and  $\widetilde{D}^{(7)}$, gives 


\begin{equation}\label{eq:zeUpdate}
   \widetilde{E}^{(7)}= \left[ \begin{array}{c} \tilde E_{(i,j,k)}^{(7)} \\  \tilde E_{(i+1,j,k)}^{(7)} \\ \vdots \\   \tilde E_{(nx,ny,nz)}^{(7)} \end{array} \right]
    = \begin{bmatrix} 
      \xi_{(i,j,k)} & 0 & \dots & 0  \\
      0 &  \xi_{(i+1,j,k)} & \dots & 0  \\
       \vdots & \vdots & \ddots & \vdots  \\
      0 & 0 & \dots & \xi_{(nx,ny,nz)}  
\end{bmatrix}
 \left[ \begin{array}{c} \tilde D_{(i,j,k)}^{(7)} \\  \tilde D_{(i+1,j,k)}^{(7)} \\ \vdots \\   \tilde D_{(nx,ny,nz)}^{(7)}  \end{array} \right],
\end{equation}

where the corresponding material matrix $\widetilde{M}_{\xi}$ is  block diagonal SPD and is ordered according to the computational cell numbering, but the global arrays consisting of the vector components for each vertex $m$,  $\widetilde{D}^{(m)}$ and $ \widetilde{E}^{(m)}$,  are ordered differently from the standard ordering. For example, for computational cell $(i,j,k)$ the standard ordering electric field components attributed to $\widetilde{E}^{(1)}$ and $\widetilde{E}^{(7)}$ respectively are, 

\begin{equation}\label{eq:zeroUpdate}
\tilde E_{(i,j,k)}^{(1)} = \left[ \begin{array}{c} E_{x(i,j,k)} \\ E_{y(i,j,k)} \\ E_{z(i,j,k)} \end{array} \right], \;\;\;\;\;   \tilde E_{(i,j,k)}^{(7)} = \left[ \begin{array}{c} E_{x(i+1,j+1,k+1)} \\ E_{y(i+1,j+1,k+1)} \\ E_{z(i+1,j+1,k+1)} \end{array} \right].
\end{equation}

Repeating this argument for each vertex labeling of the global electric field we obtain the relations

 \begin{equation} \label{eq:EkDk}
\widetilde{E}^{(m)}=\widetilde{M}_{\xi}\;\widetilde{D}^{(m)} , \;m=1,2,\dots , 8,     
\end{equation}
where the material matrix $\widetilde{M}_{\xi}$ is ordered according to the computational cell numbering which is independent of the vertex chosen to label the global electric field.  
Since the ordering of the vector components for each $\widetilde{D}^{(m)}$ and $ \widetilde{E}^{(m)}$  are distinct from the standard ordering, they  have to be permuted into a standard order denoted by  $ D$ and $E^{(m)}$, respectively,  before they can be summed up to obtain the global $ E$ field.

\smallskip
This can be accomplished by multiplying vectors $\widetilde{E}^{(m)}$ and $\widetilde{D}^{(m)}$ by corresponding orthogonal permutation matrix  $P^{(m)}$, distinct for each $m$.  Let $\displaystyle \widetilde{D}^{(m)}= P^{(m)}  D $ and $\displaystyle \widetilde{E}^{(m)}= P^{(m)} E^{(m)}$, where $D$ is a global displacement vector and $E^{(m)}$ represents the portion of the global electric field according to the formula \eqref{eq:EzAveUpdate} that  can be written in vector form as 
\begin{equation}
E =\frac{1}{8}  \sum_{m=1}^8 E^{(m)},  
\end{equation}

Substituting  $\displaystyle \widetilde{D}^{(m)}= P^{(m)}  D $ and $\displaystyle \widetilde{E}^{(m)}= P^{(m)} E^{(m)}$ into equation \ref{eq:EkDk} gives 

\begin{equation} \label{eq:permConst}
  P^{(m)} E^{(m)} = \widetilde{M}_{\xi} P^{(m)}  D,
\end{equation}  
or   
 \begin{equation} \label{eq:permConstTilde} 
 E^{(m)}=M_{\xi,m}\; D. 
\end{equation}  
where $\displaystyle M_{\xi,m} =  (P^{(m)})^T \;  \widetilde{M}_{\xi} \; P^{(m)}$. Note, that   $\displaystyle M_{\xi,m}$ matrices are  SPD matrices for each $m$, since they are orthogonally similar to an SPD matrix $\widetilde{M}_{\xi}$. 
With all of the field components now being in standard order,  the electric field update in index form in equation \ref{eq:EzAveUpdate} can be written in a vector form as

\begin{equation}
E =\frac{1}{8}  \sum_{m=1}^8 E^{(m)}=\Big(\frac{1}{8} \sum_{m=1}^8 M_{\xi,m}\Big) \;D =M_{\xi}  D,  
\end{equation}
where the global material matrix $\displaystyle M_{\xi} $ is an SPD matrix since it is equal to the sum of the SPD matrices. 
Using equation  \eqref{eq:HzAveUpdate}, a similar argument can be applied  to show that the magnetic permeability matrix $M_{\zeta}$ is a SPD matrix as well.  

\section{Numerical Validation of Stability}
\label{sec:nstability}
In this section we  validate an implementation of the fully anisotropic methods using two numerical tests. First, we perform long-time integration for a range of random high-contrast material interfaces. Second, we utilize an eigenvalue analysis of the global update matrix on a number of sample domains, showing that their eigenvalue spectra remain on the unit circle. We have not tested the diagonalizability of the global update matrix due to unstable/inconclusive nature of the Jordan form computation in finite precision arithmetic, \cite{golub13}.

\subsection{Long-Time Integration Test}

To perform the long-time integration test,  we ran simulations on a 24x24x24 cubic domain with periodic boundary conditions with lattice constant of $a = 4.8 \underline{\mu} m$ for 60 million iterations. The simulated geometries ranged from off-centered cubes and spheres to a random distribution, consisting of high-contrast, fully anisotropic materials ($\Gamma = 1, 50, 100, 144$ in equation \eqref{eq:epsAndMuBase}) embedded in a vacuum, 
\begin{equation} 
 \underline{\epsilon}(\Gamma) = \Gamma \begin{pmatrix}
    10.225 &-0.825 &-0.55 \sqrt{\frac{3}{2}} \\
    -0.825 & 10.225 & 0.55 \sqrt{\frac{3}{2}} \\
    -0.55\sqrt{\frac{3}{2}} & 0.55 \sqrt{\frac{3}{2}} & 9.95
    \end{pmatrix} \epsilon_0 , \hspace{3mm}
     \underline{\mu} (\Gamma) = \Gamma \begin{pmatrix}
    3.75 & 0.75 &-0.5 \sqrt{\frac{3}{2}} \\
    0.75 & 3.75 & -0.5 \sqrt{\frac{3}{2}} \\
    -0.5\sqrt{\frac{3}{2}} & -0.5 \sqrt{\frac{3}{2}} & 3.5
    \end{pmatrix} \mu_0.
\label{eq:epsAndMuBase}
\end{equation} 
In the random distribution each computational cell was randomly assigned as either fully anisotropic in $\underline{\epsilon}$, fully anisotropic in $\underline{\mu}$, fully anisotropic in both, or vacuum. 

We excited the domain with a broad-band,  2fs Gaussian pulse, exciting all possible modes of the empty cavity, from the lowest mode dictated by the lattice constant to the highest mode dictated by Nyquist. The computed solutions remained bounded with no evidence of exponentially growing oscillations up to the 60 millionth iteration, corresponding to approximately 2 million periods of the lowest vacuum mode.
\subsection{Eigenvalue Analysis}
Here, we demonstrate that the update matrix has eigenvalues lying on the unit circle for an arbitrary distribution of anisotropic materials. 

Both update methods being linear operators may be represented as
\begin{equation}
U^{n+1}  =  AU^n,
\end{equation}
where $A$ is the update matrix, $U$ is a vector containing the $D$ and $B$ fields at each spatial location in the computational domain, and the superscript $n$ denotes the iteration number in time. 
We numerically construct the matrix $A$ by applying a single iteration of the field update to each of the canonical basis elements of the vector $U$, e.g. $U=(1,0,\dots , 0),\; (0,1,\dots ,0)\; (0,0,\dots, , 1)$. 
The corresponding eigenvalue equation for $A$ is 
\begin{equation}
A\nu_i  =  \lambda_i \nu_i ,
\end{equation}
where $\nu_i$ is an eigenvector of $A$ and $\lambda_i$ is its corresponding eigenvalue.
For a $N\times N\times N$ domain the vector contains 6 field components from $D$ and $B$ resulting in a matrix $A$ with $(6N^3)^2$ elements, producing $6N^3$ eigenvalues.
Computing the eigenvalues of the large matrix $A$ puts a practical upper limit on the size of the computational domains on which this stability analysis can be performed. 

\begin{figure}[htb]
\begin{subfigure}{0.24\textwidth}
\includegraphics[width=1.0\linewidth]{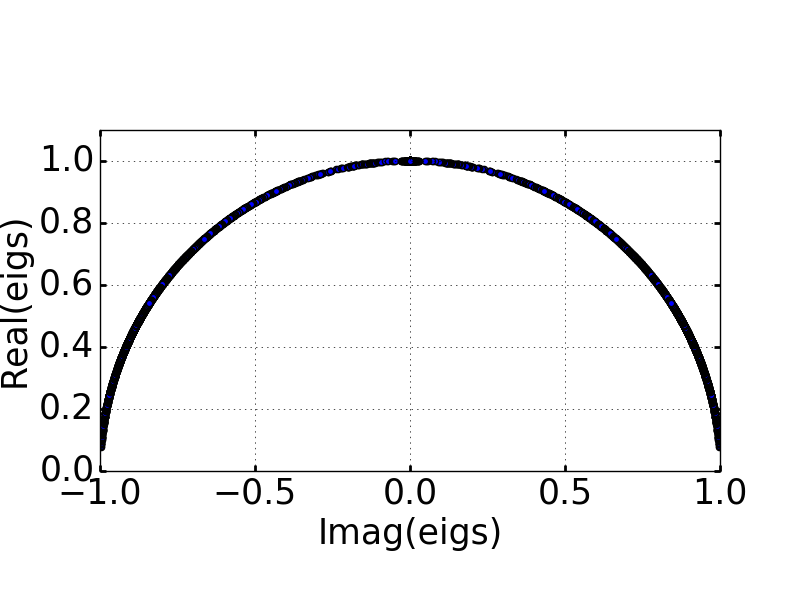}
\caption{}  \label{fig:1a}
\end{subfigure}
\hspace*{\fill} 
\begin{subfigure}{0.24\textwidth}
\includegraphics[width=1.0\linewidth]{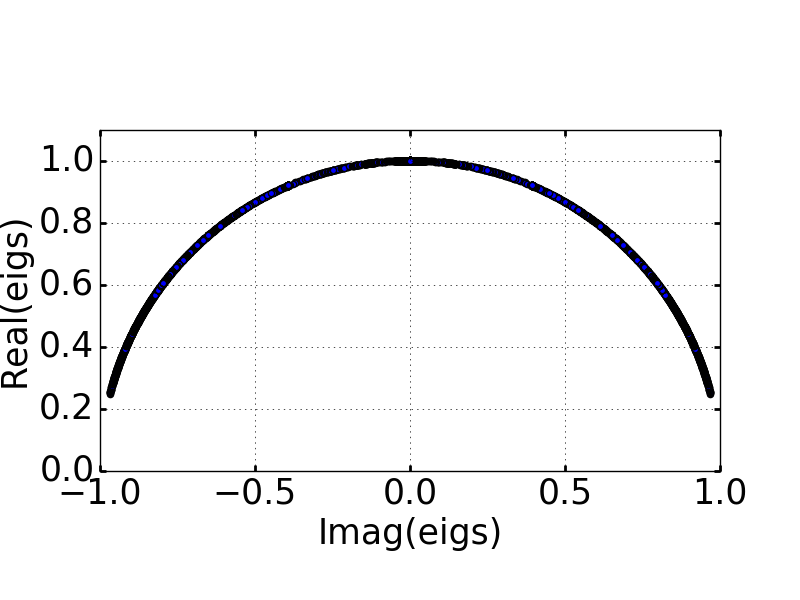}
\caption{}  \label{fig:1b}
\end{subfigure}
\begin{subfigure}{0.24\textwidth}
\includegraphics[width=1.0\linewidth]{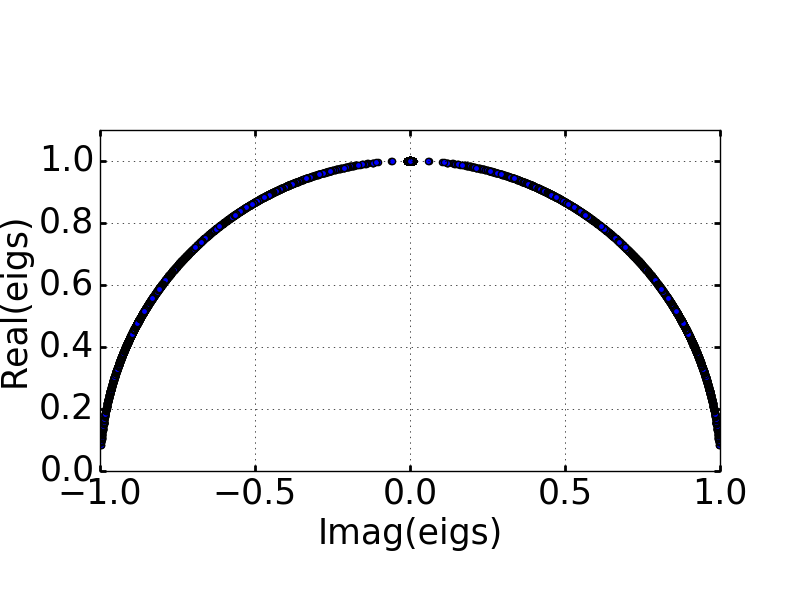}
\caption{} \label{fig:1c}
\end{subfigure}
\hspace*{\fill} 
\begin{subfigure}{0.24\textwidth}
\includegraphics[width=1.0\linewidth]{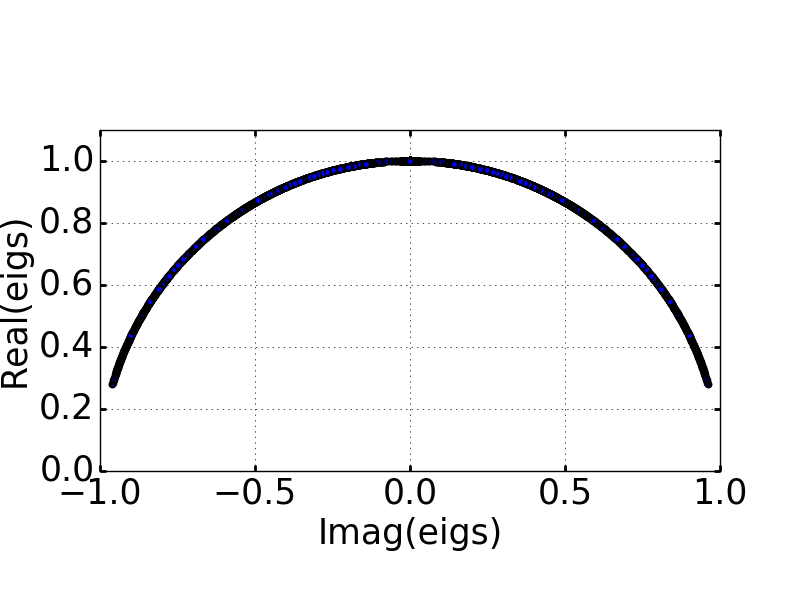}
\caption{}  \label{fig:1d}
\end{subfigure}
\caption{Eigenvalues of the updated matrix plotted for various material layouts of the high-contrast fully anisotropic materials:  (a) non-averaged method with an asymmetric sphere,   (b) non-averaged method with a random distribution of materials, (c) averaged method with an asymmetric sphere,   (d) averaged method with a random distribution of materials. In each test (including the geometries not shown here), all of the eigenvalues remained on the unit circle. }
\label{fig:eigenPlots}
\end{figure}

The update matrix, $A$, encodes information regarding the entire computational domain configuration, including material properties and layout, CFL, and boundary conditions; thus the eigenvalue stability analysis is only relevant to the computational domains tested. We numerically construct update matrices and compute their eigenvalues for computational domains with a range of materials and material distributions. In the following, we present the results for two representative material layouts: an isolated asymmetrically located sphere and a random distribution of high-contrast fully anisotropic materials. 
Using high-contrast material interfaces in test domains is particularly important due to their tendency to introduce instability \cite{werner2013more}.
Smooth material distributions, such as those used as metamaterial cloaks, were also tested with the same, stable results. 
A sphere is chosen because, when discretized, it presents a complex staircased material interface with a wide range of nearest neighbor material arrangements.
Similarly, a domain composed of a random distribution of small, high-index structures presents a large material interface surface, maximizing the possibility of seeding exponentially growing solutions. 
For the analysis, $12\times12\times12$ grids are used, resulting in $10368\times10368$ update matrices, and the CFL used was $0.4$.
The eigenvalues corresponding to each of the update matrices are calculated and plotted along with the unit circle in figure \ref{fig:eigenPlots}.
In each case we note that each of the $10368$ eigenvalues of $A$ lie on the unit circle within machine precision of the eigenvalue solver.

\section{Accuracy and Convergence}
\label{sec:accuracy}

Having proved the stability and validated the implementation numerically of both fully anisotropic methods under consideration, we now turn our attention to their accuracy properties. Due to the limited availability of analytical test solutions for fully anisotropic materials, we construct a test inspired by the method of manufactured solutions, the key idea of which to construct an exact solution, without being concerned about its physical realism \cite{roache2002code}.
Using the transformation optics approach, we construct a complex, spatially inhomogeneous geometrical object using fully anisotropic materials such that a plane-wave  incident on the object is perfectly reconstructed downstream. Such an object was described in the field of transformation optics by Pendry \emph{et al.} \cite{pendry2006controlling} where the objects are referred to as electromagnetic cloaks.

As the test was designed for numerical investigation of the accuracy properties, we are not concerned with any physical cloaking properties of the objects other than their ability to reconstruct the plane-wave downstream. Such cloaking objects have a number of desirable properties: they are composed of fully anisotropic materials, contain materials with a large dynamic range of $\underline{\epsilon}$ and $\underline{\mu}$, and exhibit a simple plane wave solution against which we can compute  relative error. 

\begin{figure}[htb]
\begin{subfigure}{\textwidth}
  \centering
  \includegraphics[width=320px]{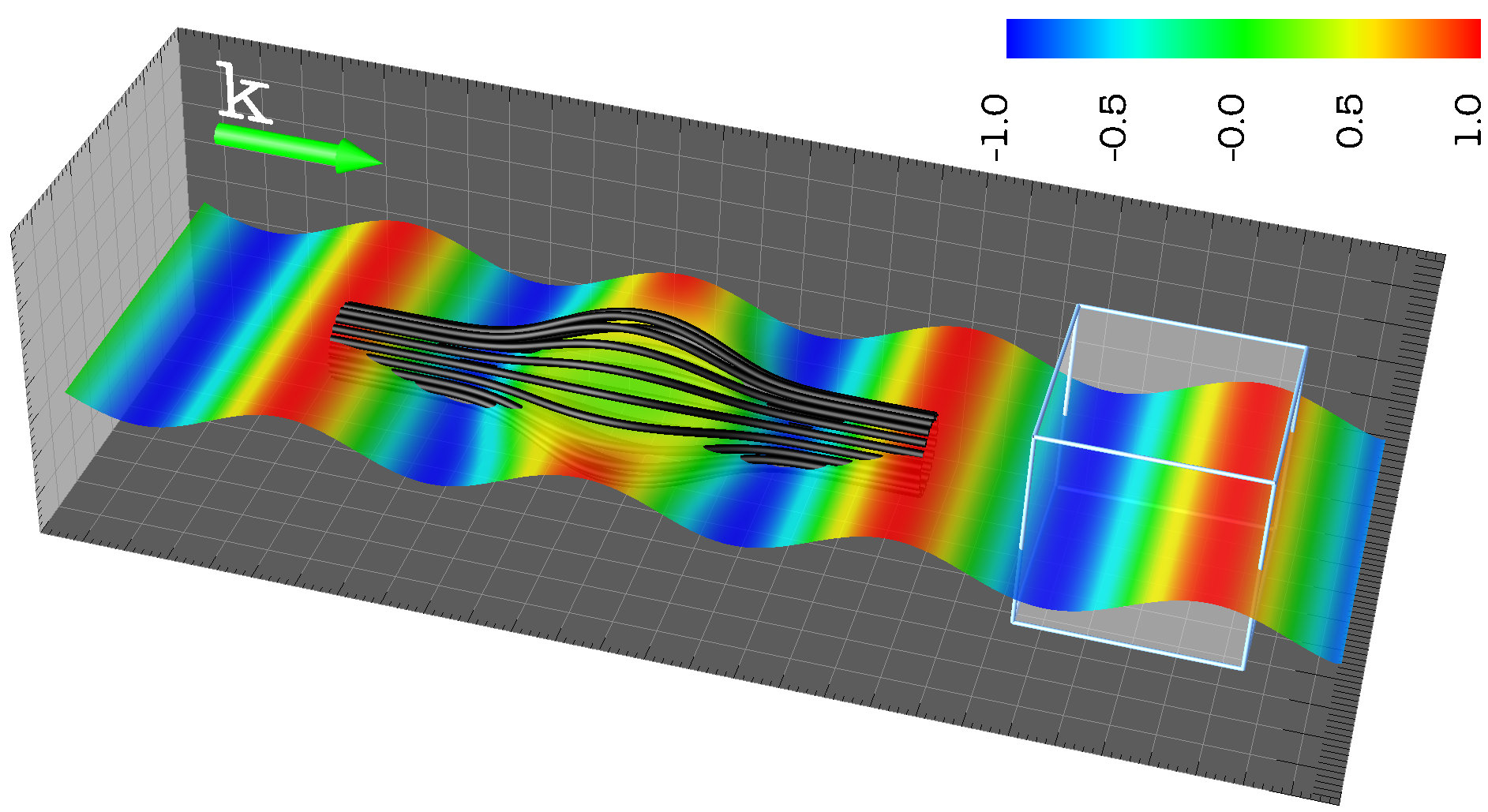}
\end{subfigure}
\caption{A visualization of the computational domain used to calculate the accuracy using a metamaterial cloak.
  The figure shows a 2D slice of the $E_y$-field, which is driven form the left side of the domain as a plane wave with amplitude 1.
  Black stream lines follow the pointing vectors of the total field in the region of the cloak, located at the center of the domain.
  The white box represents the sampling region where the relative error is computed. This figure was generated using data from of a smooth cloak using the averaged mathod at 58 ppw.} 
  \label{fig:cloakField}
\end{figure}

To characterize the field error due to discretization,
we illuminate a cloak from one side with a plane wave, and quantify the field error by comparing the computed solution to the analytically expected plane-wave over a rectangular region on the downstream side of the cloak, as shown in figure \ref{fig:cloakField}.

Using the transformation optics approach, we define two versions of a cloak: one with a smoothly-varying material distribution function with a continuous first derivative, and one for the discontinuous material distributions.

\begin{figure}[htb!]
\begin{subfigure}{0.5\textwidth}
\includegraphics[width=0.9\linewidth]{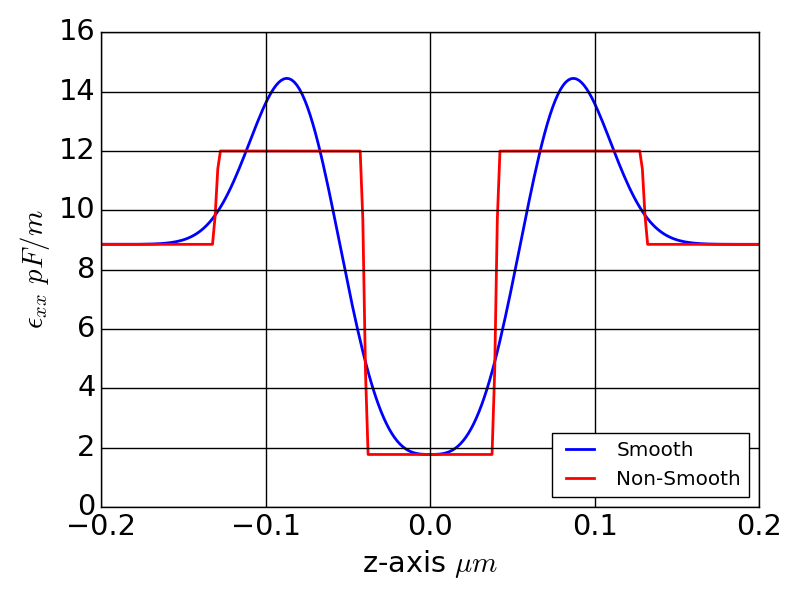}
\end{subfigure}
\hspace*{\fill} 
\begin{subfigure}{0.5\textwidth}
\includegraphics[width=0.9\linewidth]{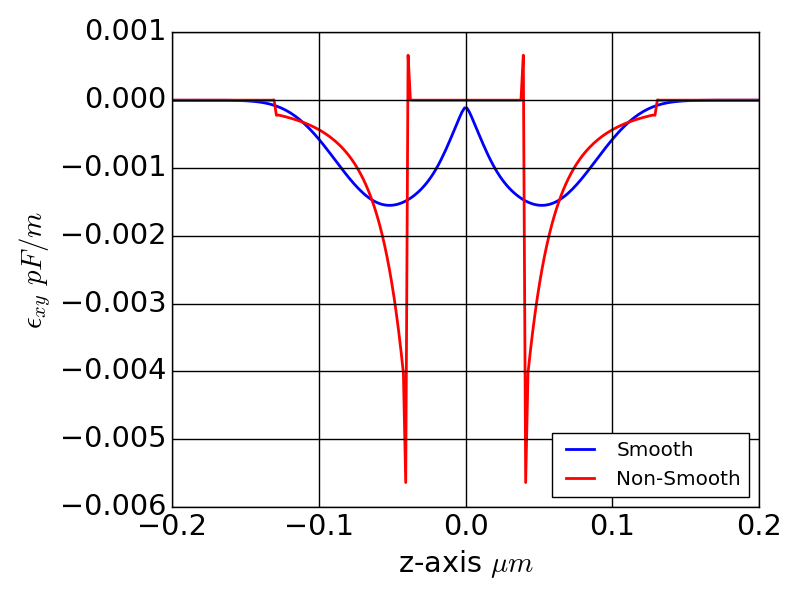}
\end{subfigure}
\caption{A 1D cut of the $\epsilon_{xx}$ (left) and $\epsilon_{xy}$ (right) along the propagation axis for the smoothly-varying (labeled Smooth) and non-smoothly-varying (labeled Non-Smooth) cloaks.}
\label{fig:epsilons}
\end{figure}

 Details of the spatial transformations used to generate the two cloaks  are shown in \ref{ap:cloak}.
The smooth cloak is used to characterize the accuracy of smoothly-varying material distributions, and the non-smooth cloak is used to characterize the accuracy of non-smoothly-varying or high-contrast material distributions.
Both of these cloaks are discretized onto the computational grid. Matrix components $\epsilon_{xx}$ and $\epsilon_{xy}$  along a 1D cut through the computational grid   are plotted in figure \ref{fig:epsilons}.

Each cloak is simulated in a 3D computational domain of size $500nm\times500nm\times1000nm$. 
The domain is discretized at various resolutions, quantified as points per wavelength (ppw),  with periodic boundary conditions employed in $x$ and $y$, and uniaxial perfectly matched layer (UPML) boundary in $z$, the propagation direction. 
A polynomial of order 3 ($m_{pml} =3$) and 10  cells of UPML ($n_{pml} = 10$) were used in computing the grading of the UPML layers.
 The  $\sigma_{max}$ and $\kappa_{max}$ are fixed to 1 and  $8 × (m_{pml} + 1)/(n_{pml}\Delta)$ respectively, where $\Delta$ is the cell size in the UPML region. 
A spherical cloaking region with an approximate radius of influence of $200nm$ is located at the center of the domain. 
Each cell in the spherical region is assigned a particular tensor value of $\underline{\epsilon}$ and $\underline{\mu}$ dictated by  either a smooth or non-smooth function defining the cloak.
A plane wave ($\lambda = 200nm$), propagating in the negative z-direction, is generated by total-field scatter-field source, located at the edge of the domain.

To ensure the results were not sensitive to the domain boundaries and/or PML thickness, several configurations were considered when determining a suitable computational domain.
The particular simulation configuration presented here has a number of advantages that minimizes the error introduced by PML boundaries.
At the source end we use a uni-directional source which limits the reflections from the source-side PML layer. 
At the transmitted end, there is near-normal incidence of the rays impinging on the PML layer.
Also, periodic boundary conditions are used on all faces normal to the propagation direction which introduce no error into the domain.

Electric fields are recorded from each computational cell in a sampling box (white region in figure \ref{fig:cloakField}) on the transmitted side of the cloaked region. The sampling box is chosen to encloses half of a wavelength. 
We then compare the electric fields in the sampling region to the analytical solution in two stages.
First,  knowing that the isotropic FDTD method converges with second-order accuracy to the analytical plane wave solution, we record a numerical simulation of a plane wave in a vacuum using the isotropic FDTD algorithm.
Second, we simulate a cloak with the fully anisotropic FDTD algorithm and compute the relative error to the  corresponding isotropic FDTD plane-wave simulation on the same grid.

We compute the scaled $L_1$ norm of the relative error using,
\begin{equation}\label{eq:relError}
  Relative Error = \frac{1}{N} \sum_j \frac{| E_j - \hat{E}_j|}{|E_j|},
\end{equation}
where $N$ is the number of points in the sampling box, $j$ indexes each point in the box, $E_j$ is the reference vacuum electric field at point $j$,  $\hat{E}_j$ is the electric field of the test run containing the cloak at point $j$.
Thus, we can determine the accuracy of the fully anisotropic FDTD to the analytical plane wave up to second-order, the order of accuracy of the isotropic FDTD to the analytical plane wave solution.

\begin{figure}[htb]
\includegraphics[width=300px]{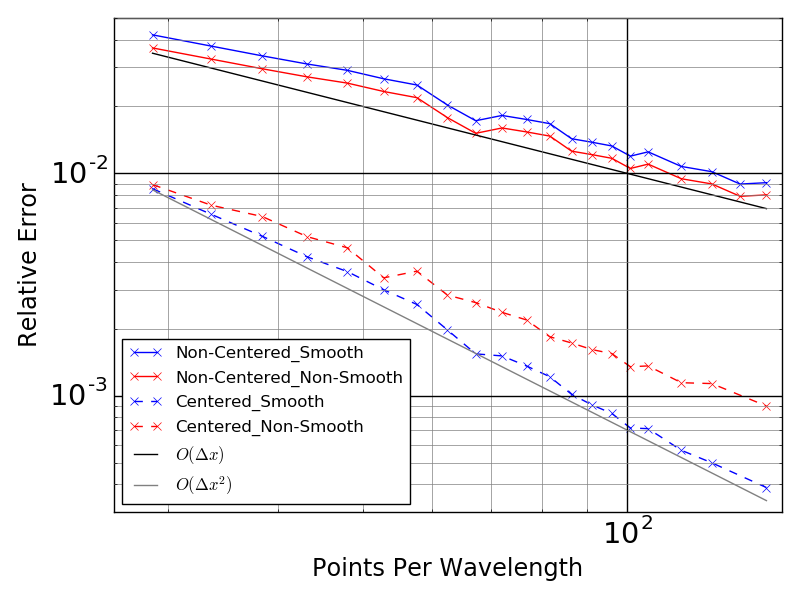}
\centering\caption{The scaled $L_1$ norm of the relative error in the sampling region for the non-averaged and averaged methods (prefixed accordingly in the legend) plotted vs ppw of the vacuum wavelength for both non-smoothly-varying and smoothly-varying cloaks (suffixed accordingly in the legend). }
\label{fig:converg}
\end{figure}

We compute the scaled $L_1$ norm of the relative error over the sampling region and plot the relative error versus the grid resolution for both the non-averaged and averaged methods in figure \ref{fig:converg}. 
The non-averaged method exhibits first-order accuracy for each cloak, as is expected due to the first-order error in material assignment on the discrete grid. 
With the averaged method, we have demonstrated second-order accuracy for the smoothly-varying cloak. 
This is consistent to what was seen in the fully anisotropic $\underline{\epsilon}$ case for continuously varying material values\cite{werner2013more}.
In our implementation, the considerable increase in accuracy and higher-order accuracy of the averaged method only increases the computational run time by 20\%.

\section{Summary}
We have developed a new FDTD Maxwell solver for fully anisotropic electric and magnetic materials assuming that the material is piecewise constant and are co-located over the primary computational cells. We have provided a mathematically rigorous proof of the stability of the two algorithms under consideration, and illustrated it by extensive numerical tests that include long-time integration, an eigenvalue analysis for extreme and random material parameter values, and various boundary conditions; accuracy evaluation for a test case designed to have an explicit analytic solution. 

We found that the non-averaged method is first order accurate, while the averaged method is second order accurate for materials with continuous first derivative and it reduces to first order for the discontinuous material distributions. In our implementation, the considerable increase in accuracy of the averaged method only increases the computational run time by 20\%.

\section{Acknowledgments}
This material is based on upon work supported by the Air Force Office of Scientific Research under award numbers FA9550-16-1-0088 and FA9550-16-1-0199.

\section*{References}

\bibliography{mybibfile}

\appendix

\section{Electromagnetic Cloak} \label{ap:cloak}
To construct the cloaking objects used in section \ref{sec:accuracy}, we follow the general principles were outlined in \cite{pendry2006controlling}. The cloaks are created by specifying a spatial transformation to 3D uniform grids that pull ray-paths out a desired region (located at the center of the domain), calculating the Jacobian matrix for the transformation, and combining the Jacobian into the material parameters ($\underline{\epsilon}$ and $\underline{\mu}$). 

For the smooth cloak, the smooth spatial transformation is defined by
\begin{equation} \label{eq:gamma}
 r^\prime = \Bigg\{1-depth\times \exp\Bigg(-\Big( \frac{r}{\sigma}\Big)^n \Bigg) \Bigg\} r,
\end{equation}
where $r$ is the uniform spatial coordinate and $r'$ is the transformed coordinate.

For the non-smooth cloak, the non-smooth piecewise-linear transformation outlined in \cite{liu2014transformation},
\begin{equation}\label{eq:lmr}
  r = f(r') =\begin{cases}
   \frac{ R_1 r'}{R'_1},  &  \text{ if } r' <R'_1, \\
   \frac{ R_1 - R_2}{R'_1- R_2}(r' - R'_1) + R_1,  &  \text{ if } R'_1 \leq r' \leq R_2,\\
   r' &  \text{ otherwise. }
\end{cases}
\end{equation}
This is used over the canonical piecewise-linear cloak \cite{pendry2006controlling} due to the infinite material parameters needed for the singularity of an ideal cloak by mapping the finite volume $R_1$ to another finite volume $R^\prime_1$ rather than a single point \cite{schurig2006metamaterial}.

We then obtain the material parameters for the cloak from the above transformations using the formulation outlined in \cite{liu2014transformation} as follows:
\begin{eqnarray}  \label{eq:transEpsAndMu}
  \underline{\epsilon}' &=& |\Lambda| \Lambda^{-1} \underline{\epsilon} \Lambda^{-T} \\
  \underline{\mu}'      &=& |\Lambda| \Lambda^{-1} \underline{\mu}      \Lambda^{-T},  
\end{eqnarray}
where  $\Lambda$ is the Jacobian matrix, which we compute numerically. 

The parameters for the  smooth cloak in the accuracy test runs were $n =3$, $depth = 0.8,$ and $\sigma = 80nm$ in equation \eqref{eq:gamma}. The parameters for the non-smooth cloak used in the accuracy test runs were $R_1 = 8nm$, $R_2 = 130nm$ and $R^\prime_1 = 40nm$ in equation \ref{eq:lmr}. Note we employ cloaks that are non-ideal in that they allow a percentage of light into the cloaked region.  This is required to keep the materials finite \cite{schurig2006metamaterial}, though under plane wave illumination, the downstream solution is the desired plane wave.

\section{Non-Averaged Constitutive Relations Update} \label{ap:non-centeredUpdate}
 
Each material and field component has a subscript $i,j,k$ that indicates the position of FDTD computational cell in the overall rectangular domain.
The relative positions of each material and field component are shown in figure \ref{fig:yeeCube}.
Furthermore, each material also has a $\{xx,xy,xz,...\}$ subscript denoting the specific matrix element. 
The $E$ and $H$ field component update are defined as follows.
\subsection{E-field update}
$E_x$:
\begin{eqnarray*}
cx &=& \xi_{xx(i,j,k)} D_{x(i,j,k)}\\ 
cy &=& \xi_{xy(i,j,k)} D_{y(i,j,k)}\\ 
cz &=& \xi_{xz(i,j,k)} D_{z(i,j,k)}\\ 
E_{x(i,j,k)} &=& cx + cy + cz
\end{eqnarray*}

$E_y$:
\begin{eqnarray*}
cx &=& \xi_{yx(i,j,k)} D_{x(i,j,k)}\\ 
cy &=& \xi_{yy(i,j,k)} D_{y(i,j,k)}\\ 
cz &=& \xi_{yz(i,j,k)} D_{z(i,j,k)}\\ 
E_{y(i,j,k)} &=& cx + cy + cz
\end{eqnarray*}

$E_z$:
\begin{eqnarray*}
cx &=& \xi_{zx(i,j,k)} D_{x(i,j,k)}\\ 
cy &=& \xi_{zy(i,j,k)} D_{y(i,j,k)}\\ 
cz &=& \xi_{zz(i,j,k)} D_{z(i,j,k)}\\ 
E_{z(i,j,k)} &=& cx + cy + cz
\end{eqnarray*}

\subsection{H-field update}
$H_x$:
\begin{eqnarray*}
cx &=& \zeta_{xx(i,j,k)} B_{x(i,j,k)}\\ 
cy &=& \zeta_{xy(i,j,k)} B_{y(i,j,k)}\\ 
cz &=& \zeta_{xz(i,j,k)} B_{z(i,j,k)}\\ 
H_{x(i,j,k)} &=& cx + cy + cz
\end{eqnarray*}

$H_y$:
\begin{eqnarray*}
cx &=& \zeta_{yx(i,j,k)} B_{x(i,j,k)}\\ 
cy &=& \zeta_{yy(i,j,k)} B_{y(i,j,k)}\\ 
cz &=& \zeta_{yz(i,j,k)} B_{z(i,j,k)}\\ 
H_{y(i,j,k)} &=& cx + cy + cz
\end{eqnarray*}

$H_z$:
\begin{eqnarray*}
cx &=& \zeta_{zx(i,j,k)} B_{x(i,j,k)}\\ 
cy &=& \zeta_{zy(i,j,k)} B_{y(i,j,k)}\\ 
cz &=& \zeta_{zz(i,j,k)} B_{z(i,j,k)}\\ 
H_{z(i,j,k)} &=& cx + cy + cz
\end{eqnarray*}
\section{Averaged Constitutive Relations Update} \label{ap:aveUpdate}
In this section we explicitly describe field and material averaging methods for fully anisotropic $\underline{\epsilon}$ and $\underline{\mu}$. 

\subsection{E-field update}\label{ap:aveUpdateE}
$E_x$:
\begin{eqnarray*}
cx &=& \frac{1}{2}(\xi_{xx(i,j,k)}+\xi_{xx(i-1,j,k)}) Dx_{x(i,j,k)}\\ 
cy &=& \frac{1}{4}(\xi_{xy(i,j,k)} (D_{y(i,j,k)} + D_{y(i,j+1,k)}) + \\ 
   &&        \xi_{xy(i-1,j,k)} (D_{y(i-1,j,k)} + D_{y(i-1,j+1,k)}))\\ 
cz &=& \frac{1}{4}(\xi_{xz(i,j,k)} (D_{z(i,j,k)} + D_{z(i,j,k+1)}) +\\  
   &&        \xi_{xz(i-1,j,k)} (D_{z(i-1,j,k)} + D_{z(i-1,j,k+1)}))\\ 
E_{x(i,j,k)} &=& cx + cy + cz
\end{eqnarray*}

$E_y$:
\begin{eqnarray*}
cx &=& \frac{1}{4}(\xi_{yx(i,j,k)} (D_{x(i,j,k)} + D_{x(i+1,j,k)}) + \\ 
   &&        \xi_{yx(i,j-1,k)} (D_{x(i,j-1,k)} + D_{x(i+1,j-1,k)}))\\ 
cy &=& \frac{1}{2}(\xi_{yy(i,j,k)} + \xi_{yy(i,j-1,k)}) D_{y(i,j,k)}\\ 
cz &=& \frac{1}{4}(\xi_{yz(i,j,k)} (D_{z(i,j,k)} + D_{z(i,j,k+1)}) + \\ 
   &&        \xi_{yz(i,j-1,k)} (D_{z(i,j-1,k)} + D_{z(i,j-1,k+1)})) \\ 
E_{y(i,j,k)} &=& cx + cy + cz
\end{eqnarray*}

$E_z$:
\begin{eqnarray*}
cx &=& \frac{1}{4}(\xi_{zx(i,j,k)} (D_{x(i,j,k)} + D_{x(i+1,j,k)}) + \\ 
   &&        \xi_{zx(i,j,k-1)} (D_{x(i,j,k-1)} + D_{x(i+1,k,k-1)}))\\ 
cy &=& \frac{1}{4}(\xi_{zy(i,j,k)} (D_{y(i,j,k)} + D_{y(i,j+1,k)}) +\\  
   &&        \xi_{zy(i,j,k-1)} (D_{y(i,j,k-1)} + D_{y(i,j+1,k-1)}))\\ 
cz &=& \frac{1}{2}(\xi_{zz(i,j,k)} + \xi_{zz(i,j,k-1)}) D_{z(i,j,k)}\\ 
E_{z(i,j,k)} &=& cx + cy + cz
\end{eqnarray*}

\subsection{H-field update}\label{ap:aveUpdateH}
$H_x$:
\begin{eqnarray*}
cx &=& \frac{1}{4}(\zeta_{xx(i,j,k)}+\zeta_{xx(i,j-1,k)}+\zeta_{xx(i,j,k-1)}+\zeta_{xx(i,j-1,k-1)}) B_{x(i,j,k)}\\ 
cy &=& \frac{1}{8}((\zeta_{xy(i,j,k)} + \zeta_{xy(i,j,k-1)}) (B_{y(i,j,k)} + B_{y(i+1,j,k)}) +\\ 
   &&      (\zeta_{xy(i,j-1,k)} + \zeta_{xy(i,j-1,k-1)}) (B_{y(i,j-1,k)} + B_{y(i+1,j-1,k)}))\\ 
cz &=& \frac{1}{8}((\zeta_{xz(i,j,k)} + \zeta_{xz(i,j-1,k)}  ) (B_{z(i,j,k)} + B_{z(i+1,j,k)}) +\\ 
   &&         (\zeta_{xz(i,j,k-1)} + \zeta_{xz(i,j-1,k-1)}) (B_{z(i,j,k-1)} + B_{z(i+1,k,k-1)}))\\ 
H_{x(i,j,k)} &=& cx + cy + cz
\end{eqnarray*}

$H_y$:
\begin{eqnarray*}
cx &=& \frac{1}{8}((\zeta_{yx(i,j,k)} + \zeta_{yx(i,j,k-1)}) (B_{x(i,j,k)} + B_{x(i,j+1,k)}) +\\ 
   &&         (\zeta_{yx(i-1,j,k)} + \zeta_{yx(i-1,j,k-1)}) (B_{x(i-1,j,k)} + B_{x(i-1,j+1,k)}))\\ 
cy &=& \frac{1}{4}(\zeta_{yy(i,j,k)} + \zeta_{yy(i-1,j,k)}+\zeta_{yy(i,j,k-1)} + \zeta_{yy(i-1,j,k-1)}) B_{y(i,j,k)}\\ 
cz &=& \frac{1}{8}((\zeta_{yz(i,j,k)} + \zeta_{yz(i-1,j,k)}) (B_{z(i,j,k)} + B_{z(i,j+1,k)}) +\\ 
   &&         (\zeta_{yz(i,j,k-1)} + \zeta_{yz(i-1,j,k-1)}) (B_{z(i,j,k-1)} + B_{z(i,j+1,k-1)}))\\ 
H_{y(i,j,k)} &=& cx + cy + cz
\end{eqnarray*}

$H_z$:
\begin{eqnarray*}
cx &=& \frac{1}{8}((\zeta_{zx(i,j,k)}+\zeta_{zx(i,j-1,k)}) (B_{x(i,j,k)} + B_{x(i,j,k+1)}) +\\ 
   &&         (\zeta_{zx(i-1,j,k)} + \zeta_{zx(i-1,j-1,k)}) (B_{x(i-1,j,k)} + B_{x(i-1,j,k+1)}))\\ 
cy &=& \frac{1}{8}((\zeta_{zy(i,j,k)}+\zeta_{zy(i-1,j,k)}) (B_{y(i,j,k)} + B_{y(i,j,k+1)}) +\\ 
   &&         (\zeta_{zy(i,j-1,k)}+\zeta_{zy(i-1,j-1,k)}) (B_{y(i,j-1,k)} + B_{y(i,j-1,k+1)}))\\ 
cz &=& \frac{1}{4}(\zeta_{zz(i,j,k)} + \zeta_{zz(i-1,j,k)}+\zeta_{zz(i,j-1,k)} + \zeta_{zz(i-1,j-1,k)}) B_{z(i,j,k)}\\ 
H_{z(i,j,k)} &=& cx + cy + cz
\end{eqnarray*}

\end{document}